\documentclass[11pt]{amsart}
\usepackage{inputenc}
\usepackage{amsmath,amssymb, amsthm,delarray,amsopn,amscd,mathrsfs}
\usepackage{graphicx}
\numberwithin{equation}{section}
\newtheorem{thm}{Theorem}[section]

\newtheorem{definition}{Definition}{\rm}
\newtheorem{example}{Example}{\rm}
{\rm}

\def\beq{\begin{equation} }
\def\eeq{\end{equation} }

\def\N{\mathbb{N}}
\def\M{\mathbf{M}}
\def\T{\mathbf{T}}
\def\R{\mathbb{R}}
\def\A{\mathbf{A}}
\def\P{\mathbf{P}}
\def\bL{\mathbf{L}}

\def\K{\mathbf{K}}
\def\B{\mathbf{B}}
\def\Q{\mathbf{Q}}
\def\M{\mathbf{M}}

\def\X{\mathbf{X}}
\def\U{\mathbf{U}}
\def\G{\mathbf{G}}

\def\v{\mathbf{v}}
\def\f{\mathbf{f}}

\def\u{\mathbf{u}}
\def\f{\mathbf{f}}
\def\b{\mathbf{b}}

\def\x{\mathbf{x}}

\def\y{\mathbf{y}}
\def\z{\mathbf{z}}
\def\t{\mathbf{t}}
\def\s{\mathcal{S}}
\def\bS{\mathbf{S}}

\def\om{\mathbf{\Omega}}

\def\R{\mathbb{R}}
\def\B{\mathbf{B}}
\def\C{\mathbf{C}}
\def\X{\mathbf{X}}

\def\supmu{{\rm supp}(\mu)}
\begin{document}
\title[optimization]{The Moment-SOS Hierarchy}

\author{Jean B. Lasserre}

\thanks{Research supported by the European Research Council (ERC) through ERC-Advanced Grant \# 666981 for the TAMING project}

\begin{abstract}
The Moment-SOS hierarchy initially introduced in optimization in 2000, is based on the 
theory of the $\K$-moment problem and its dual counterpart, polynomials that are positive on $\K$.
It turns out that this methodology can be also applied to solve 
problems with positivity constraints ``$f(\x)\geq0$ for all $\x\in\K$"
and/or linear constraints on Borel measures. Such problems can be viewed as specific instances of
the ``Generalized Problem of Moments" (GPM) whose list of important applications in various domains
is endless. We describe this methodology and outline some of its applications in various domains.
\end{abstract}


\maketitle


\section{Introduction}
~

Consider the optimization problem:
\begin{equation}
\label{def-pb}
\P:\quad f^*\,=\,\displaystyle\inf_\x\,\{\,f(\x):\: \x\in\om\,\},
\end{equation}
where $f$ is a polynomial and $\om\subset\R^n$ is a basic semi-algebraic set, that is,
\begin{equation}
\label{setk}
\om\,:=\,\{\,\x\in\R^n:\quad g_j(\x)\,\geq\,0,\quad j=1,\ldots,m\,\},
\end{equation}
for some polynomials $g_j$, $j=1,\ldots,m$.  
Problem $\P$ is a particular case of {\em Non Linear Programming} (NLP) where the data ($f,g_j$, $j=1,\ldots,m$) 
are {\em algebraic}, and therefore the whole arsenal of methods of NLP can be used for solving $\P$. 
So what is so specific about $\P$ in (\ref{def-pb})? 
The answer depends on the meaning of $f^*$ in (\ref{def-pb}). 

If one is interested in a {\em local minimum} only then
efficient NLP methods can be used for solving
$\P$. In such methods, the fact that $f$ and $g_j$'s are polynomials does not help much, that is,
this algebraic feature of $\P$ is not really exploited.
On the other hand if $f^*$ in (\ref{def-pb}) is understood as the {\em global minimum} of $\P$ then the picture is totally different.
Why? First, to eliminate any ambiguity on the meaning of $f^*$ in (\ref{def-pb}), rewrite (\ref{def-pb}) as:
\begin{equation}
\label{def-pb-1}
\P:\quad f^*\,=\,\sup\,\{\,\lambda:\: f(\x)-\lambda\,\geq\,0,\quad \forall \x\in\om\,\}
\end{equation}
because then indeed $f^*$ is necessarily the global minimum of $\P$. 

In full generality, most problems (\ref{def-pb-1}) are very difficult to solve (they are labelled NP-hard in the computational complexity terminology) because:
\begin{center}
{\em Given $\lambda\in\R$, checking whether
``$f(\x)-\lambda\geq0$ for all $\x\in\om$" is difficult.}
\end{center}
Indeed, by nature this positivity constraint is {\em global} and therefore cannot be 
handled by standard NLP optimization algorithms which use only 
local information around a current iterate $\x\in\om$. Therefore to compute $f^*$ in (\ref{def-pb-1}) one needs an efficient tool to handle the positivity constraint 
``$f(\x)-\lambda\geq0$ for all $\x\in\om$". 
Fortunately if the data are algebraic then:
\begin{enumerate}
\item Powerful {\em positivity certificates} from Real Algebraic Geometry
 ({\em Posi-tivstellens\"atze} in german) are available.
\item Some of these positivity certificates have an efficient practical implementation via {\em Linear Programming} (LP)
or {\em Semidefinite Programming} (SDP). In particular and importantly, testing whether a given polynomial is a sum of squares (SOS) simply reduces to solving  a single SDP (which  can be done in time polynomial in the input size of the polynomial, up to arbitrary fixed precision).
\end{enumerate}
After the pioneers works of Shor \cite{shor}
and Nesterov \cite{nesterov}, Lasserre \cite{lass-cr,lass-siopt} 
and Parrilo  \cite{parrilo-1,parrilo-2} 
have been the first to provide a systematic use of 
these two key ingredients in Control and Optimization, with convergence guarantees. 
It is also worth mentioning another closely related pioneer work, namely the celebrated SDP-relaxation 
of Goemans \& Williamson \cite{goemans} which provides a  $0.878$ approximation guarantee for MAXCUT,
a famous problem in non-convex combinatorial optimization (and probably the simplest one). In fact it is perhaps the first
famous example of such a successful application of the powerful SDP convex optimization technique 
to provide guaranteed good approximations to a notoriously difficult non-convex optimization problem. It turns out that this SDP relaxation is the first relaxation in the 
Moment-SOS hierarchy (a.k.a. Lasserre hierarchy) when applied to the MAXCUT problem.
Since then, this spectacular success story of SDP relaxations has been at the origin of a 
flourishing  research activity in
combinatorial optimization and computational complexity. In particular,
the study of LP- and SDP-relaxations in hardness of approximation is at the core of
a central topic in combinatorial optimization and computational complexity, namely proving/disproving Khot's famous Unique Games Conjecture\footnote{For this conjecture and its theoretical and practical implications, S. Khot was awarded the prestigious Nevanlinna prize at the last ICM 2014 in Seoul \cite{khot2} .}  (UGC) in Theoretical Computer Science.

Finally, another ``definition" of the global optimum $f^*$ of $\P$ reads:
\begin{equation}
\label{def-pb-2}
f^*\,=\,\inf_\mu\,\{\,\int_\om f\,d\mu:\quad \mu(\om)\,=\,1\,\}
\end{equation}
where the `$\inf$" is over all probability measures on $\om$. 
Equivalently, writing $f$ as $\sum_\alpha f_\alpha\,\x^\alpha$ in the basis of monomials (where $\x^\alpha=x_1^{\alpha_1}\cdots x_n^{\alpha_n}$):
\begin{equation}
\label{def-pb-3}
f^*\,=\,\inf_\y\,\{\,\sum_\alpha f_\alpha \,y_\alpha: \y \in \mathscr{M}(\om);\quad y_0=1\,\},
\end{equation}
where $\mathscr{M}(\om)=\{\y=(y_\alpha)_{\alpha\in\N^n}:\:\exists \,\mu \mbox{ s.t.  } y_\alpha=\int_\om \x^\alpha\,d\mu,\:\forall\alpha\in\N^n\}$, a convex cone. In fact
(\ref{def-pb-1}) is the LP dual of (\ref{def-pb-2}). In other words 
standard LP duality between the two formulations (\ref{def-pb-2}) and (\ref{def-pb-1}) 
illustrates the duality between the ``$\om$-moment problem" and ``polynomials positive on $\om$".

Problem (\ref{def-pb-2}) is a very particular instance (and even the simplest instance) of the more general {\em Generalized Problem of Moments} (GPM):
\begin{equation}
\label{GPM-primal}
\inf_{\mu_1,\ldots,\mu_p}\{\,\sum_{j=1}^p\int_{\om_j} f_j\,d\mu_j:\:
\sum_{j=1}^p f_{ij}\,d\mu_j\,\geq\, b_i,\:i=1,\ldots,s\,\},
\end{equation}
for some functions $f_{ij}:\R^{n_j}\to\R$, $i=1,\ldots,s$, and sets $\om_j\subset\R^{n_j}$, $j=1,\ldots,p$. The GPM is an infinite-dimensional LP with dual:
\begin{equation}
\label{GPM-dual}
\sup_{\lambda_1,\ldots,\lambda_s\ge0}\{\,\sum_{i=1}^s \lambda_i\,b_i:\:
f_j-\sum_{i=1}^s \lambda_{i}\,f_{ij}\,\geq\,0\mbox { on $\om_j$},\:j:1,\ldots,p\}.
\end{equation}
Therefore it should be of no surprise that the Moment-SOS hierarchy, initially developed for global optimization, also applies to solving the GPM. This is particularly interesting as the list of important applications of the GPM is almost endless; see e.g. Landau \cite{moments}.

\section{The MOMENT-SOS hierarchy in optimization}

\subsection{Notation, definitions and preliminaries}

Let $\R[\x]$ denote the ring of polynomials in the variables $\x=(x_1,\ldots,x_n)$ and let
$\R[\x]_d$ be the vector space of polynomials of degree at most $d$
(whose dimension is $s(d):={n+d\choose n}$).
For every $d\in\N$, let  $\N^n_d:=\{\alpha\in\N^n:\vert\alpha\vert \,(=\sum_i\alpha_i)\leq d\}$, 
and
let $\v_d(\x)=(\x^\alpha)$, $\alpha\in\N^n$, be the vector of monomials of the canonical basis 
$(\x^\alpha)$ of $\R[\x]_{d}$. 
Given a closed set $\mathcal{X}\subseteq\R^n$, let $\mathscr{P}(\mathcal{X})\subset\R[\x]$ (resp. $\mathscr{P}_d(\mathcal{X})\subset\R[\x]_d$) be the convex cone of polynomials (resp. polynomials of degree at most $2d$) that are nonnegative on $\mathcal{X}$. A polynomial $f\in\R[\x]_d$ is written
\[\x\mapsto f(\x)\,=\,\sum_{\alpha\in\N^n}f_\alpha\,\x^\alpha,\]
with vector of coefficients $\f=(f_\alpha)\in\R^{s(d)}$ in the canonical basis of monomials $(\x^\alpha)_{\alpha\in\N^n}$.
For real symmetric matrices, let $\langle \B,\C\rangle:={\rm trace}\,(\B\C)$ while the notation $\B\succeq0$ 
stands for $\B$ is positive semidefinite (psd) whereas 
$\B\succ0$ stands for $\B$ is  positive definite (pd).

\subsection*{The Riesz functional}
Given a sequence $\y=(y_\alpha)_{\alpha\in\N^n}$, the Riesz functional is the linear mapping
$L_\y:\R[\x]\to\R$ defined by:
\begin{equation}
\label{Riesz}
f\:(=\sum_\alpha f_\alpha\,\x^\alpha)\quad \mapsto L_\y(f)\,=\,\sum_{\alpha\in\N^n}f_\alpha\,y_\alpha.
\end{equation}

\subsection*{Moment matrix}
The {\it moment} matrix associated with a sequence
$\y=(y_\alpha)$, $\alpha\in\N^n$, is the real symmetric matrix $\M_d(\y)$ with rows and columns indexed by $\N^n_d$, and whose entry 
$(\alpha,\beta)$ is just $y_{\alpha+\beta}$, for every $\alpha,\beta\in\N^n_d$. 
Alternatively, let
$\v_d(\x)\in\R^{s(d)}$ be the vector $(\x^\alpha)$, $\alpha\in\N^n_d$, and
define the matrices $(\B_{o,\alpha})\subset\s^{s(d)}$ by
\begin{equation}
\label{balpha}
\v_d(\x)\,\v_d(\x)^T\,=\,\sum_{\alpha\in\N^n_{2d}}\B_{o,\alpha}\,\x^\alpha,\qquad\forall\x\in\R^n.\end{equation}
Then $\M_d(\y)=\sum_{\alpha\in\N^n_{2d}}y_\alpha\,\B_{o,\alpha}$.
If $\y$ has a representing measure $\mu$ then
$\M_d(\y)\succeq0$ because $\langle\f,\M_d(\y)\f\rangle=\int f^2d\mu\geq0$, for all $f\in\R[\x]_d$.

A measure whose all moments are finite, is {\it moment determinate} if there is no other measure with same moments.
The support of a Borel measure $\mu$ on $\R^n$ (denoted $\supmu$) is the smallest closed set $\om$ such that $\mu(\R^n\setminus\om)=0$.

\subsection*{Localizing matrix}
With $\y$ as above and $g\in\R[\x]$ (with $g(\x)=\sum_\gamma g_\gamma\x^\gamma$), the {\it localizing} matrix associated with $\y$ 
and $g$ is the real symmetric matrix $\M_d(g\,\y)$ with rows and columns indexed by $\N^n_d$, and whose entry $(\alpha,\beta)$ is just $\sum_{\gamma}g_\gamma y_{\alpha+\beta+\gamma}$, for every $\alpha,\beta\in\N^n_d$.
Alternatively, let $\B_{g,\alpha}\in\s^{s(d)}$ be defined by:
\begin{equation}
\label{calpha}
g(\x)\,\v_d(\x)\,\v_d(\x)^T\,=\,\sum_{\alpha\in\N^n_{2d+{\rm deg}\,g}}\B_{g,\alpha}\,\x^\alpha,\qquad\forall\x\in\R^n.\end{equation}
Then $\M_d(g\,\y)=\sum_{\alpha\in\N^n_{2d+{\rm deg}g}}y_\alpha\,\B_{g,\alpha}$.
If $\y$ has a representing measure $\mu$ whose support is 
contained in the set $\{\x:g(\x)\geq0\}$ then
$\M_d(g\,\y)\succeq0$ for all $d$ because $\langle\f,\M_d(g\,\y)\f\rangle=\int f^2\,gd\mu\geq0$,
for all $f\in\R[\x]_d$.

\subsection*{SOS polynomials and quadratic modules}
A polynomial $f\in\R[\x]$ is a Sum-of-Squares (SOS) if  there exist $(f_k)_{k=1,\ldots,s}\subset\R[\x]$, such that
$f(\x)=\sum_{k=1}^s f_k(\x)^2$, for all $\x\in\R^n$. Denote by $\Sigma[\x]$ (resp. $\Sigma[\x]_d$) the set of 
SOS polynomials (resp. SOS polynomials of degree at most $2d$).
Of course every SOS polynomial is nonnegative whereas the converse is not true. 
In addition, checking  whether 
a given polynomial $f$ is nonnegative on $\R^n$ is difficult whereas checking whether $f$ is SOS is much easier
and can be done efficiently. Indeed let $f\in\R[\x]_{2d}$ (for $f$ to be SOS its degree must be even), $\x\mapsto f(\x)=\sum_{\alpha\in\N^n_{2d}}f_\alpha\,\x^\alpha$. Then $f$ is SOS if and only if there exists a real  symmetric matrix $\X^T=\X$ of size $s(d)={n+d\choose n}$, such that:
\begin{equation}
\label{sos}
\X\succeq0;\quad f_\alpha\,=\,\langle \X,\B_{o,\alpha}\rangle,\qquad \forall \alpha\in\N^n_{2d},
\end{equation}
and this can be checked by solving an SDP.

Next, let $\x\mapsto g_0(\x):=1$ for all $\x\in\R^n$. With a family $(g_1,\ldots,g_m)\subset\R[\x]$ is associated 
the {\it quadratic module} $Q(g)\,(=Q(g_1,\ldots,g_m))\,\subset\R[\x]$:
\begin{equation}
\label{module}
Q(g)\,:=\,\left\{\sum_{j=0}^m\sigma_j\,g_j:\:\sigma_j\in\Sigma[\x],\:j=0,\ldots,m\,\right\},
\end{equation}
and its {\it truncated} version
\begin{equation}
\label{truncated-module}
Q_k(g)\,:=\,\left\{\sum_{j=0}^m\sigma_j\,g_j:\: \sigma_j\in\Sigma[\x]_{k-d_j},\:j=0,\ldots,m\,\right\},
\end{equation}
where $d_j=\lceil{\rm deg}(g_j)/2\rceil$, $j=0,\ldots,m$.

\begin{definition}
The quadratic module $Q(g)$ associated with $\om$ in (\ref{setk}) is said to be Archimedean if there exists $M>0$ such that
the quadratic polynomial $\x\mapsto M-\Vert\x\Vert^2$ belongs to  $Q(g)$ (i.e., belongs to $Q_k(g)$ for some $k$).
\end{definition}
If $Q(g)$ is Archimedean then necessarily  $\om$ is compact but the reverse is not rue. 
The Archimedean condition (which depends on the representation of $\om$) can be seen as an {\em algebraic certificate} that $\om$ is compact.
For more details on the above notions of moment and localizing matrix, quadratic module,
as well as their use in potential applications, the interested reader is referred to Lasserre \cite{lass-book-icp}, Laurent \cite{laurent-2},
Schm\"udgen \cite{schmudgen}.

\subsection{Two certificates of positivity (Positivstellens\"atze)}

Below we describe two particular certificates of positivity which are important because they provide the theoretical justification
behind the so-called SDP- and LP-relaxations for global optimization.
\begin{thm}[Putinar \cite{putinar}]
\label{th-put}
Let $\om\subset\R^n$ be as in (\ref{setk}) and assume that $Q(g)$ is Archimedean. 

(a) If a  polynomial $f\in\R[\x]$ is (strictly) positive on $\om$ then $f\in Q(g)$.

(b) A sequence $\y=(y_\alpha)_{\alpha\in\N^n}\subset\R$ has a representing Borel measure on $\om$
if and only if $L_\y( f^2\,g_j)\geq0$ for all $f\in\R[\x]$, and all $j=0,\ldots,m$. Equivalently, if and only if $\M_d(\y\,g_j)\succeq0$ for all $j=0,\ldots,m$, $d\in\N$.
\end{thm}
There exists another certificate of positivity which does not use SOS. 
\begin{thm}[Krivine-Vasilescu \cite{krivine1,krivine2,vasilescu}]
\label{th-krivine}
Let $\om\subset\R^n$ as in (\ref{setk}) be compact and such that (possibly after scaling)
$0\leq g_j(\x)\leq 1$ for all $\x\in\om$, $j=1,\ldots,m$. Assume also that $[1,g_1,\ldots,g_m]$ generates $\R[\x]$.
 
 (a) If a  polynomial $f\in\R[\x]$ is (strictly) positive on $\om$ then 
 \begin{equation}
 \label{th-krivine-1}
 f(\x)\,=\,\sum_{\alpha,\beta\in\N^n}c_{\alpha,\beta}\,\prod_{j=1}^m g_j(\x)^{\alpha_j}\,(1-g_j(\x))^{\beta_j},
 \end{equation}
 for finitely many positive coefficients $(c_{\alpha,\beta})_{\alpha,\beta\in\N^m}$.
 
 (b) A sequence $\y=(y_\alpha)_{\alpha\in\N^n}\subset\R$ has a representing Borel measure on $\om$ if and only if
 $L_\y\left(\displaystyle\prod_{j=1}^m g_j(\x)^{\alpha_j}\,(1-g_j(\x))^{\beta_j}\right)\,\geq\,0$ for all $\alpha,\beta\in\N^m$.
 \end{thm}
The two facets (a) and (b) of Theorem \ref{th-put} and Theorem \ref{th-krivine} illustrate the duality between {\em polynomials
positive on $\om$} (in (a)) and the {\em $\om$-moment problem} (in (b)).
In addition to their mathematical interest, both Theorem 2.1(a) and Theorem \ref{th-krivine}(a)  have another distinguishing feature. They both 
have a practical implementation.
Testing whether $f\in\R[\x]_d$ is in $Q(g)_k$ is just solving a single SDP, whereas testing whether $f$ van be written as in
(\ref{th-krivine-1}) with $\sum_{i=1}^m\alpha_i+\beta_i\leq k$, is just solving a single Linear Program (LP).

\subsection{The Moment-SOS hierarchy}

The Moment-SOS hierarchy is a numerical scheme based on Putinar's theorem. In a nutshell it consists of replacing the intractable positivity constraint
``$f(\x)\geq0$ for all $\x\in\om$" with Putinar's positivity certificate $f\in Q_d(g)$ of Theorem \ref{th-put}(a), i.e., with a fixed degree bound on the SOS weights $(\sigma_j)$ in (\ref{truncated-module}). By duality, it consists of replacing the intractable constraint $\y\in\mathscr{M}(\om)$ with the necessary conditions 
$\M_d(g_j\,\y)\succeq0$, $j=0,\ldots,m$, of Theorem \ref{th-put}(b) for  a fixed $d$.
This results in solving an SDP which provides a lower bound on the global minimum. By allowing the degree bound $d$ to increase,
one obtains a {\it hierarchy}  of SDPs (of increasing size)  which provides a monotone non-decreasing sequence of lower bounds.
A similar strategy based on Krivine-Stengle-Vasilescu positivity certificate (\ref{th-krivine-1}) is also possible and yields a hierarchy of LP (instead of SDPs). However
even though one would prefer to solve LPs rather than SDPs, the latter Moment-LP hierarchy has several serious drawbacks 
(some explained in e.g. \cite{lass-book-camb,lass-mor}), and therefore we only describe the Moment-SOS hierarchy.

Recall problem $\P$ in (\ref{def-pb}) or equivalently in (\ref{def-pb-1}) and (\ref{def-pb-2}), where $\om\subset\R^n$ is the basic semi-algebraic set
defined in (\ref{setk}). 

\subsection*{The Moment-SOS hierarchy} Consider the sequence of semidefinite programs $(\Q_d)_{d\in\N}$ 
with $d\geq \hat{d}:=\max[{\rm deg}(f), \max_j {\rm deg}(g_j)]$:
\begin{equation}
\label{sdp-moment}
\Q_d:\: \rho_d=\displaystyle\inf_\y \,\{\,L_\y(f):\:y_0=1;\:\M_d(g_j\,\y)\,\succeq0,\quad 0\leq j\leq m\,\}
\end{equation}
(where $\y=(y_\alpha)_{\alpha\in\N^n_{2d}}$)\footnote{In Theoretical Computer Science,  $\y$ is called a sequence of ``pseudo-moments".}, with associated sequence of  their SDP duals:
\begin{equation}
\label{sdp-sos}
\Q^*_d:\:
\rho^*_d=\displaystyle\sup_{\lambda,\sigma_j} 
\{\,\lambda:\:f-\lambda =\displaystyle\sum_{j=0}^m\sigma_j\,g_j;\: \sigma_j\in\Sigma[\x]_{d-d_j},\: 0\leq j\leq m\}
\end{equation}
(where $d_j=\rceil ({\rm deg} g_j)/2\rceil$). 
By standard weak duality in optimization $\rho_d^*\leq\rho_d$ for every $d\geq \hat{d}$.
The sequence $(\Q_d)_{d\in\N}$ forms a {\em hierarchy} of {\em SDP-relaxations} of $\P$ because $\rho_d\leq f^*$ and $\rho_d\leq \rho_{d+1}$ for all $d\geq \hat{d}$.
Indeed for each $d\geq\hat{d}$,  the constraints of $\Q_d$ consider only necessary conditions for $\y$ to be the moment sequence (up to order $2d$) of a probability measure on $\om$ (cf. Theorem \ref{th-put}(b)) and therefore $\Q_d$ is a relaxation of (\ref{def-pb-3}). 

By duality, the sequence $(\Q^*_d)_{d\in\N}$ forms a {\em hierarchy} of {\em SDP-strenghtenings} of (\ref{def-pb-1}). Indeed in (\ref{sdp-sos}) 
one has  replaced the intractable positivity constraint of (\ref{def-pb-1}) by the (stronger) Putinar's positivity certificate with degree bound $2d-2d_j$ on the SOS weights $\sigma_j$'s.
\begin{thm}[\cite{lass-cr,lass-siopt}]
Let $\om$ in (\ref{setk}) be compact and assume that its associated quadratic module $Q(g)$ is Archimedean. Then:

(i) As $d\to\infty$, the monotone non-decreasing 
sequence $(\rho_d)_{d\in\N}$ (resp. $(\rho^*_d)_{d\in\N}$) of optimal values of the hierarchy (\ref{sdp-moment}) 
(resp. (\ref{sdp-sos})) converges to the global optimum $f^*$ of $\P$.

(ii) Moreover, let $\y^d=(y^d_\alpha)_{\alpha\in\N^n_{2d}}$ be an optimal solution of $\Q_d$ in (\ref{sdp-moment}), and let $s=\max_jd_j$ (recall that  
$d_j=\lceil( {\rm deg}\,g_j)/2\rceil$).  If
\begin{equation}
\label{rank}
{\rm rank}\,\M_d(\y^d)\,=\,{\rm rank}\,\M_{d-s}(\y^d)\:(=:t)
\end{equation}
then $\rho_d=f^*$ and there are $t$ global minimizers $\x^*_j\in\om$, $j=1,\ldots,t$, that can be ``extracted" from $\y^d$
by a linear algebra routine.
\end{thm}
The sequence of SDP-relaxations $(\Q_d)$, $d\geq\hat{d}$, and the rank test (\ref{rank}) to extract global minimizers, 
are implemented in the GloptiPoly software \cite{gloptipoly}.

\subsection*{Finite convergence and a global optimality certificate}

After being introduced in \cite{lass-cr}, 
in many numerical experiments it was observed that typically, finite convergence takes place, that is, $f^*=\rho_d$ for some 
(usually small) $d$. In fact there is a rationale behind this empirical observation. 
\begin{thm}[Nie \cite{nie-1}]
\label{th-nie-global}
Let $\P$ be as in (\ref{def-pb-1}) where $\om$ in (\ref{setk}) is compact and its associated quadratic module is Archimedean. Suppose that at each global minimizer $\x^*\in\om$:

$\bullet$ The gradients $(\nabla g_j(\x^*))_{j=1,\ldots,m}$ are linearly independent. (This implies existence of nonnegative 
Lagrange-KKT multipliers $\lambda^*_j$, $j\leq m$, such that $\nabla f(\x^*)-\sum_{j=1}^m\lambda_j^*\,\nabla g_j(\x^*)=0$
and $\lambda^*_jg_j(\x^*)=0$ for all $j\leq m$.)

$\bullet$ Strict complementarity holds, that is, $g_j(\x^*)=0\Rightarrow \lambda_j^*>0$.

$\bullet$ Second-order sufficiency condition holds, i.e., \\
\[\langle\u,\nabla_\x^2 \,(f(\x^*)-\sum_{j=1}^m\lambda_j^*\,g_j(\x^*))\,\u\rangle\,>\,0,\]
for all $0\neq\u\in \nabla (f(\x^*)-\sum_{j=1}^m\lambda_j^*\,g_j(\x^*))^\perp$.

Then $f-f^*\in Q(g)$, i.e., there exists $d^*$ and SOS multipliers $\sigma^*_j\in\Sigma[\x]_{d^*-d_j}$, $j=0,\ldots,m$, such that:
\begin{equation}
\label{th-nie-global-1}
f(\x)-f^*\,=\,\sigma^*_0(\x)+\sum_{j=1}^m \sigma^*_j(\x)\,g_j(\x).
\end{equation}
\end{thm}
With (\ref{th-nie-global-1}), Theorem \ref{th-nie-global} provides a 
{\em certificate of global optimality} in polynomial optimization, and to the best of our knowledge, the first at this level of generality.
Next, observe that $\x^*\in\om$ is a global unconstrained minimizer of the {\em extended Lagrangian polynomial} $f-f^*-\sum_{j=1}^n\sigma^*_jg_j$,
and therefore Theorem \ref{th-nie-global} is the  analogue
 for {\em non-convex} polynomial optimization of the Karush-Kuhn-Tucker (KKT) optimality conditions {\em in the convex case}. Indeed in the convex case, any local minimizer is global
 and is also a global unconstrained minimizer of the Lagrangian $f-f^*-\sum_{j=1}^m\lambda^*_j g_j$. 
 
 Also interestingly, whenever the SOS weight $\sigma^*_j$ in (\ref{th-nie-global-1}) is non trivial, it testifies that the constraint $g_j(\x)\geq0$ is important for $\P$
 even if it is not active  at $\x^*$ (meaning that if $g_j\geq0$ is deleted from $\P$ then the new global optimum decreases  strictly).  The multiplier $\lambda^*_j$ plays
 the same role in the KKT-optimality conditions {\em only in the convex case}. See \cite{lass-book-camb} for a detailed discussion.
 
\noindent
{\bf Finite convergence} of the Moment-SOS-hierarchies (\ref{sdp-moment}) and (\ref{sdp-sos}) is an immediate consequence of Theorem 
\ref{th-nie-global}. Indeed by (\ref{th-nie-global-1}) $(f^*,\sigma_0^*,\ldots,\sigma_m^*)$ is a feasible solution of $\Q^*_{d^*}$ with value
$f^*\leq\rho^*_d\leq f^*$ (hence $\rho^*_d=\rho_d=f^*$).

\noindent
{\bf Genericity}: Importantly,  as proved in Nie \cite{nie-1}, the conditions in 
Theorem \ref{th-nie-global} are {\em generic}. By this we mean the following: Consider the class $\mathscr{P}(t,m)$ of optimization problems $\P$
with data $(f,g_1,\ldots,g_m)$ of degree bounded by $t$, and with nonempty compact feasible set $\om$. Such a problem $\P$ is a ``point"
in the space $\R^{(m+1)s(t)}$ of coordinates of $(f,g_1,\ldots,g_m)$. Then the ``good" problems $\P$ are points 
in a Zariski open set.
Moreover, generically the rank test (\ref{rank}) is also satisfied at an optimal solution of (\ref{sdp-moment}) (for some $d$); for more details see Nie \cite{nie-2}.

\noindent
{\bf Computational complexity}:
Each relaxation $\Q_d$ in (\ref{sdp-moment}) is a semidefinite program with $s(2d)={n+2d\choose n}$ variables $(y_\alpha)$, and 
a psd constraint $\M_d(\y)\succeq0$ of size $s(d)$. Therefore solving $\Q_d$ in its canonical form (\ref{sdp-moment})
is quite expensive in terms of computational burden, especially when using interior-point methods. Therefore 
its brute force application is limited to small to medium size problems. 

\noindent
{\bf Exploiting sparsity}: Fortunately many 
large scale problems exhibit a structured sparsity pattern (e.g., each polynomial $g_j$ is concerned with a few variables only, and the objective function $f$ is a sum $\sum_if_i$ where each $f_i$ 
is also concerned with a few variables only). Then Waki et al. \cite{waki} have proposed
a sparsity-adapted hierarchy of SDP-relaxations which can handle problems $\P$ with thousands variables. In addition, if the sparsity pattern satisfies a certain condition then convergence of this sparsity-adapted hierarchy is also guaranteed like in the dense case \cite{lass-sparsity}. Successful applications of this strategy can be found in e.g. Camps and Sznaier \cite{sznaier} in Control (systems identification) and in Molzahn and Hiskens \cite{molzahn} for solving (large scale) 
Optimum Power Flow problems (OPF is an important problem encountered in the management of energy networks).

\subsection{Discussion}

We claim that the Moment-SOS hierarchy and its rationale Theorem \ref{th-nie-global}, unify convex, non-convex (continuous), and discrete (polynomial)  Optimization.
Indeed in the description of $\P$ we do not pay attention to what particular  class of  problems $\P$ belongs to. This is in sharp contrast to the usual 
common practice in (local) optimization where several classes of problems have their own tailored favorite class of algorithms.
For instance, problems are not treated the same if equality constraints appear, and/or if boolean (or discrete variables)  are present, etc. Here a 
boolean variable $x_i$ is modeled by the quadratic equality constraint $x_i^2=x_i$. So it is reasonable to speculate 
that this lack of specialization could be  a handicap for the moment-SOS hierarchy.

But this is not so. For instance for the sub-class of convex\footnote{Convex problems $\P$ where $f$ and $(-g_j)_{j=1,\ldots,m}$ are convex, are considered ``easy" and can be solved efficiently.} problems
 $\P$ where $f$ and $(-g_j)_{j=1,\ldots,m}$ are SOS-convex\footnote{A polynomial $f\in\R[\x]$ is SOS-convex if its Hessian $\nabla^2f$ is a SOS matrix-polynomial, i.e., 
$\nabla f^2(\x)=\bL(\x)\bL(\x)^T$ for some matrix-polynomial $\bL\in\R[\x]^{n\times p}$.} polynomials, finite convergence takes place at the first step of the hierarchy. In other words, the SOS hierarchy somehow ``recognizes" this class of easy problems \cite{lass-book-camb}.
In the same time, for  a large class of $0/1$ combinatorial optimization problems on graphs, the Moment-SOS hierarchy has been shown to provide the tightest upper bounds when compared to the class of {\em lift-and-project} methods, and  has now become a central tool to analyze hardness of approximations in combinatorial optimization.
For more details the interested reader is referred to e.g. Lasserre \cite{lass-mor}, Laurent \cite{laurent-1}, Barak \cite{barak}, Khot \cite{khot,khot2}
and the many references therein. 

\section{The Moment-SOS hierarchy outside optimization}

\subsection{A general framework for the Moment-SOS hierarchy}
Let $\om_i\subset\R^{n_i}$ be a finite family of compact sets, $\mathscr{M}(\om_i)$ (resp. $\mathscr{C}(\om_i)$) be the space of finite Borel signed measures (resp. continuous functions) on $\om_i$, $i=0,1,\ldots,s$, and let $\T$ be a continuous linear mapping with adjoint $\T^*$:
\begin{eqnarray*}
\T:\mathscr{M}(\om_1)\times\cdots\times \mathscr{M}(\om_s)&\to&\mathscr{M}(\om_0)\\
\mathscr{C}(\om_1)\times\cdots\times\mathscr{C}(\om_s)&\leftarrow&\mathscr{C}(\om_0):\:\T^*
\end{eqnarray*}
Let $\phi:=(\phi_1,\ldots,\phi_s)$ and let $\phi_i\geq0$ stand for $\phi_i$ is a positive measure. Then consider the general framework:
\begin{equation}
\label{frame-1}
\rho=\displaystyle\inf_{\phi\geq0}\,\{\,\displaystyle\sum_{i=1}^s \langle f_i,\phi_i\rangle:
\T(\phi)\,=\lambda; \:
\displaystyle\sum_{i=1}^s\langle f_{ij},\phi_i\rangle\,\geq\,b_j,\: j\in J\},
\end{equation}
where $J$ is a finite or countable set, $\b=(b_j)$ is given, $\lambda\in\mathscr{M}(\om_0)$ is a given measure, 
$(f_{ij})_{j\in J}$, $i=1,\ldots,s$, are given polynomials, and $\langle\cdot,\cdot\rangle$ is the duality bracket between $\mathscr{C}(\om_i)$ and $\mathscr{M}(\om_i)$ ($\langle h,\phi_i\rangle=\int_{\om_i} hd\phi_i$), $i=1,\ldots,s$. 

As we will see, this general framework is quite rich as it encompasses a lot of important applications in many different fields. In fact Problem (\ref{frame-1}) is equivalent to the Generalized Problem of Moments (GPM):
\begin{equation}
\label{frame-2}
\begin{array}{rl}
\rho=\displaystyle\inf_{\phi\geq0}&\{\,\displaystyle\sum_{i=1}^s \langle f_i,\phi_i\rangle: 
\langle \T^*\,p_k,\phi\rangle\,=\,\langle p_k,\lambda\rangle,\quad k=0,1,\ldots\\
&\displaystyle\sum_{i=1}^s\langle f_{ij},\phi_i\rangle\,\geq\,b_j,\quad j\in J\},
\end{array}
\end{equation}
where the family $(p_k)_{k=0,\ldots}$ is dense in $\mathscr{C}(\om_0)$ (e.g.
a basis of $\R[x_1,\ldots,x_{n_0}]$). 

The Moment-SOS hierarchy can also be applied to
help solve the Generalized Problem of Moments (GPM)  (\ref{frame-2}) 
or its dual :
\begin{equation}
\label{frame-2-dual}
\begin{array}{rl}
\rho^*=\displaystyle\sup_{(\theta_j\geq0,\mathbf{\gamma})}&\{\,\displaystyle\sum_k\gamma_k\,\langle p_k,\lambda\rangle+\langle \theta,\b\rangle:\\
\mbox{s.t.}&f_i-\displaystyle\sum_{k}\gamma_k\,(\T^*\,p_k)_i -\sum_{j\in J}\theta_j\,f_{ij}\geq\,0\:\mbox{ on $\om_i$ for all $i$\,}\},
\end{array}
\end{equation}
where the unknown $\gamma=(\gamma_k)_{k\in\N}$ is a finite sequence.

\subsection{A hierarchy of SDP-relaxations}
Let 
\begin{equation}
\label{omega}
\om_i\,:=\,\{\,\x\in\R^{n_i}: \:g_{i,\ell}(\x)\geq0,\,i=1,\ldots,m_i\,\},\quad i=1,\ldots,s,
\end{equation}
for some polynomials $(g_{i,\ell})\subset\R[x_{1},\ldots,x_{n_i}]$, $\ell=1,\ldots,m_i$. Let
$d_{i,\ell}=\lceil{\rm deg}(g_{i,\ell})/2\rceil$ and 
$\hat{d}:=\max_{i,j,\ell}[{\rm deg}(f_i),{\rm deg}(f_{ij}),{\rm deg}(g_{i,\ell})]$.
To solve (\ref{frame-2}), define the ``moment" sequences $\y_i=(y_{i,\alpha})$, $\alpha\in\N^{n_i}$, $i=1,\ldots,s$, and with $d\in\N$, define $\Gamma_d:=\{p_k: {\rm deg}(T^*p_k)_i\leq 2d,\:i=1,\ldots,s\}$.
Consider the hierarchy of semidefinite programs indexed by $\hat{d}\leq d\in\N$:
\begin{equation}
\label{frame-2-sdp}
\begin{array}{rl}
\rho_d=\displaystyle\inf_{(\y_i)}&\{\,\displaystyle\sum_{i=1}^s L_{\y_i}(f_i): \:\displaystyle\sum_{i=1}^s L_{\y_i}((T^*p_k)_i)\,=\,\langle p_k,\lambda\rangle,\quad p_k\in\Gamma_d\\
&\displaystyle\sum_{i=1}^sL_{\y_i}(f_{ij})\,\geq\,b_j,\quad j\in J_d\\
&\M_d(\y_i),\:\M_{d-d_\ell}(g_{i\ell}\,\y_i)\succeq0,\quad \ell\leq m_i;\:i\leq s\},
\end{array}
\end{equation}
where $J_d\subset J$ is finite $\bigcup_{d\in\N} J_d=J$. Its dual SDP-hierarchy reads:
\begin{equation}
\label{frame-2-sdp-dual}
\begin{array}{rl}
\rho^*_d=\displaystyle\sup_{(\theta_j\geq0,\gamma_k)}&\{\,\displaystyle\sum_{p_k\in\Gamma_d}\gamma_k\,\langle p_k,\lambda\rangle+\langle \theta,\b\rangle:\\
\mbox{s.t.}&f_i-\displaystyle\sum_{p_k\in\Gamma_d}\gamma_k\,(\T^*\,p_k)_i -\sum_{j\in J}\theta_j\,f_{ij}=\sum_{\ell=0}^{m_i}\sigma_{i,\ell}\,g_{i,\ell}\\
&\sigma_{i,\ell}\in\Sigma[x_{1},\ldots,x_{n_i}]_{d-d_{i,\ell}};\:
i=1,\ldots,s\},
\end{array}
\end{equation}

As  each $\om_i$ is compact, for technical reasons and with no loss of generality, in the sequel we may and will assume that for every $i=1,\ldots,s$, $g_{i,0}(\x)=M_i-\Vert \x\Vert^2$,  where $M_i>0$ is sufficiently large.

\begin{thm}
\label{th-general}
Assume that $\rho>-\infty$ and that for every $i=1,\ldots,s$, $f_{i0}=1$.
Then for every $d\geq \hat{d}$, (\ref{frame-2-sdp}) has an optimal solution, and $\lim_{d\to\infty}\rho_d=\rho$.
\end{thm}

\subsection{Examples in Probability and Computational Geometry}

\subsection*{Bounds on measures with moment conditions}
Let $Z$ be a random vector with values in a compact semi-algebraic set 
$\om_1\subset\R^n$. Its distribution 
$\lambda$ on $\om_1$ is unknown but some
of its moments 
$\int \x^\alpha\,d\lambda=b_\alpha$,  $\alpha\in\Gamma\subset\N^n$, are known ($b_0=1$).
 Given a basic semi-algebraic set $\om_2\subset\om_1$ we want to compute (or approximate as closely as desired) 
the best upper bound on ${\rm Prob}(Z\in \om_2)$. This problem reduces to solving
the GPM:
\begin{equation}
\label{ex1}
\begin{array}{rl}
\rho=\displaystyle\sup_{\phi_1,\phi_2\geq0}&\{\langle 1,\phi_2\rangle: \langle \x^\alpha,\phi_1\rangle+
\langle \x^\alpha,\phi_2\rangle=b_\alpha,\:\alpha\in\Gamma;\\
&\phi_i\in\mathscr{M}(\om_i),\:i=1,2\,\},\end{array}
\end{equation}
With $\om_1$ and $\om_2$  as in (\ref{omega}) one may compute upper bounds on $\rho$ by solving the Moment-SOS hierarchy (\ref{frame-2-sdp}) adapted to problem (\ref{ex1}). Under the assumptions of Theorem \ref{th-general},
the resulting sequence $(\rho_d)_{d\in\N}$ converges to $\rho$ as $d\to\infty$;
for more details the interested reader is referred to \cite{aap}.

\subsection*{Lebesgue \& Gaussian measures of semi-algebraic sets}

Let $\om_2\subset\R^n$ be compact. The goal is to compute (or approximate as closely as desired) the Lebesgue measure $\lambda(\om_2)$ of $\om_2$. Then take 
$\om_1\supset\om_2$ be a simple set, e.g. an ellipsoid or a box (in fact any set such that
one knows all moments $(b_\alpha)_{\alpha\in\N^n}$ of the Lebesgue measure on $\om_1$). Then:
\begin{equation}
\label{ex2}
\begin{array}{rl}
\lambda(\om_2)\,=\,\displaystyle\sup_{\phi_1,\phi_2\geq0}&\{\langle 1,\phi_2\rangle: \langle \x^\alpha,\phi_1\rangle+
\langle \x^\alpha,\phi_2\rangle=b_\alpha,\:\alpha\in\N^n;\\
&\phi_i\in\mathscr{M}(\om_i),\:i=1,2\,\}.\end{array}
\end{equation}
Problem (\ref{ex2}) is very similar to (\ref{ex1}) except that we now have countably many moment constraints ($\Gamma=\N^n$). Again, with 
$\om_2$ and $\om_2$ as in (\ref{omega}) one may compute upper bounds on 
$\lambda(\om_2)$ by solving the Moment-SOS hierarchy (\ref{frame-2-sdp}) 
adapted to problem (\ref{ex2}). Under the assumptions of Theorem \ref{th-general},
the resulting monotone non-increasing sequence $(\rho_d)_{d\in\N}$ converges to $\lambda(\om_2)$ from above as $d\to\infty$.
The convergence $\rho_d\to\lambda(\om_2)$ is slow because of a Gibb's phenomenon\footnote{The Gibbs' phenomenon appears at a jump discontinuity when one approximates a piecewise $C^1$ function with a continuous function, e.g., by its Fourier series.}. Indeed the semidefinite program (\ref{frame-2-sdp-dual}) reads:
\[\rho^*_d=\inf_{p\in\R[\x]_{2d}}\,\{ \int_{\om_1}p\,d\lambda:\:p\geq1\mbox{ on $\om_2$};\quad
p\geq0\mbox{ on $\om_1$}\},\]
i.e., as $\to\infty$ one tries to approximate the discontinuous function $\x\mapsto 1_{\om_2}(\x)$ by polynomials of increasing degrees. Fortunately there are several ways to accelerate the convergence, e.g. as in \cite{sirev} (but loosing the monotonicity) or in \cite{aam} (preserving monotonicity)  by including in (\ref{frame-2-sdp}) additional constraints on $\y_2$
coming from an application of Stokes' theorem.

For the {\bf Gaussian measure} $\lambda$ we need and may take $\om_1=\R^n$
and $\om_2$ is not necessarily compact. Although both $\om_1$ and $\om_2$ are allowed to be non-compact, the Moment-SOS hierarchy (\ref{frame-2-sdp}) still converges, i.e., $\rho_d\to\lambda(\om_2)$ as $d\to\infty$. This is because
the moments of $\lambda$ satisfy the generalized Carleman's condition
\begin{equation}
\label{carleman}
\sum_{k=1}^\infty \left(\int_{\R^n} x_i^{2k}\,d\lambda\right)^{-1/2k}\,=\,+\infty,\quad i=1,\ldots,n,
\end{equation}
which imposes implicit constraints on $\y_1$ and $\y_2$ in (\ref{frame-2-sdp}),
strong enough to guarantee 
$\rho_d\to\lambda(\om_2)$ as $d\to\infty$. For more details see \cite{aam}. 
This deterministic approach is computationally demanding and should be seen as complementary to brute force Monte-Carlo 
methods that provide only an estimate (but can handle larger size problems).

\subsection{In signal processing and interpolation}

In this application, a signal is identified 
with an atomic signed measure $\phi$ supported on few atoms $(\x_k)_{k=1,\ldots,s}\subset\om$, i.e., $\phi=\sum_{k=1}^s \theta_k\,\delta_{\x_k}$, for some weights $(\theta_k)_{k=1,\ldots,s}$. 

\subsection*{Super-Resolution} The goal of Super-Resolution is to reconstruct the unknown measure $\phi$ (the signal) from a few measurements only, when those measurements are the moments $(b_\alpha)_{\alpha\in\N^n_t}$
of $\phi$, up to order $t$ (fixed). One way to proceed is to solve the infinite-dimensional program:
\begin{equation}
\label{resolution}
\rho\,=\,\displaystyle\inf_{\phi}\,\{\Vert\phi\Vert_{TV}: \int \x^\alpha\,d\phi\,=\,b_\alpha,\quad \alpha\in\N^n_t\,\},
\end{equation}
where the $\inf$ is over the finite signed Borel measures on $\om$,
and $\Vert\phi\Vert_{TV}=\vert\phi\vert(\om)$ (with $\vert\phi\vert$ 
being the total variation of $\phi$). Equivalently:
\begin{equation}
\label{resolution-1}
\rho\,=\,\displaystyle\inf_{\phi^+,\phi^-\geq0}\,\{\langle 1,\phi^++\phi^-\rangle: \langle \x^\alpha, \phi^+-\phi^-\rangle\,=\,b_\alpha,\: \alpha\in\N^n_t\,\},
\end{equation}
which is an instance of the GPM with dual:
\begin{equation}
\label{resolution-2}
\rho^*\,=\,\displaystyle\sup_{p\in\R[\x]_t}\,\{
\sum_{\alpha\in\N^n_t}p_\alpha\,b_\alpha:
\Vert p\Vert_\infty\,\leq\,1\,\},
\end{equation}
where $\Vert p\Vert_\infty=\sup\{\vert p(\x)\vert:\x\in\om\}$. In this case, the Moment-SOS hierarchy (\ref{frame-2-sdp}) with $d\geq \hat{d}:=\lceil t/2\rceil$, reads:
\begin{equation}
\label{resolution-3}
\begin{array}{rl}.
\rho_d=\displaystyle\inf_{\y^+,\y^-}&\{ y_0^++y^-_0: 
y^+_\alpha-y^-_\alpha\,=\,b_\alpha,\: \alpha\in\N^n_t\\
&\M_d(\y^\pm)\succeq0; \:\M_d(g_\ell\,\y^\pm)\succeq0,\:\ell=1,\ldots,m\},
\end{array}
\end{equation}
where $\om=\{\x: g_\ell(\x)\geq0,\:\ell=1,\ldots,m\}$.

In the case where $\om$ is the torus $\mathbb{T}\subset\mathbb{C}$, Cand\`es and Fernandez-Granda \cite{candes} showed 
that if $\delta>2/f_c$ (where $\delta$ is the minimal distance between the atoms of $\phi$, and $f_c$ is the number of measurements) then (\ref{resolution}) has a unique solution and one may recover $\phi$ {\em exactly} by solving the single semidefinite program (\ref{resolution}) with $d=\lceil t/2\rceil$. The dual (\ref{resolution-2}) has an optimal solution $p^*$ (a trigonometric polynomial) and  the support of $\phi^+$ (resp. $\phi^-$) consists of
the atoms $\z\in\mathbb{T}$ of $\phi$ such that $p^*(\z)=1$ (resp. $p^*(\z)=-1$). In addition, this procedure is 
more robust to noise in the measurements than Prony's method; on the other hand, the latter requires less measurements and no separation condition on the atoms. 

In the general multivariate case treated in \cite{super-resolution} 
one now needs to solve the  Moment-SOS hierarchy (\ref{resolution-1}) for $d=\hat{d},\ldots$
(instead of a single SDP in the univariate case).
However since the moment constraints of (\ref{resolution-1}) are finitely many, exact
recovery (i.e. finite convergence of the Moment-SOS hierarchy (\ref{resolution-3})) is possible (usually with a few measurements only). This is indeed what has been observed in all  numerical experiments of \cite{super-resolution}, and in all cases with significantly less 
measurements than the theoretical bound 
(of a tensorized version of the univariate case). 

In fact, the rank condition (\ref{rank}) is always satisfied at an optimal solution
$(\y^+,\y^-)$ at some step $d$ of the hierarchy (\ref{resolution-3}), and
so the atoms of $\phi^+$
and $\phi^-$ are extracted via a simple linear algebra routine
(as for  global optimization). 
Nie's genericity result \cite{nie-2} should provide a rationale
which explains why the rank condition (\ref{rank}) is satisfied in all examples.

\subsection*{Sparse interpolation} Here the goal is to recover an unknown (black-box) polynomial $p\in\R[\x]_t$ through a few evaluations of $p$ only. In \cite{interpolation} we have shown that this problem is in fact a particular case of Super-Resolution (and even {\em discrete} Super-Resolution) on the torus $\mathbb{T}^n\subset\mathbb{C}^n$. Indeed let $\z_0\in\mathbb{T}^n$ be fixed, arbitrary. Then with $\beta\in\N^n$, notice that
\begin{eqnarray*}
p(\z_0^\beta)&=&\sum_{\alpha\in\N^n_d}p_\alpha\,
(z_{01}^{\beta_1}\cdots z_{0n}^{\beta_n})^\alpha
\,=\,\sum_{\alpha\in\N^n_d}p_\alpha\,
(z_{01}^{\alpha_1}\cdots z_{0n}^{\alpha_n})^\beta\\
&=&\int_{\mathbb{T}^n}\z^\beta\,d\left(\sum_{\alpha\in\N^n_d}p_\alpha\,\delta_{\z_0^\alpha}\right)\,=\,
\int_{\mathbb{T}^n} \z^\beta\,d\phi.
\end{eqnarray*}
In other words, one may identify the polynomial $p$ with an atomic  signed
Borel measure $\phi$ on $\mathbb{T}^n$ supported on finitely many atoms 
$(\z_0^\alpha)_{\alpha\in\N^n_t}$ with associated weights
$(p_\alpha)_{\alpha\in\N^n_t}$. 

Therefore, if the evaluations of the black-box polynomial $p$ are done at a few ``powers" $(\z_0^\beta)$, $\beta\in\N^n$, 
of an arbitrary point $\z_0\in\mathbb{T}^n$, then
the sparse interpolation problem is equivalent to recovering 
an unknown atomic signed Borel measure $\phi$ on $\mathbb{T}^n$
from knowledge of a few moments, that is, the Super-Resolution problem that we have just described above. Hence one may recover
$p$ by solving the Moment-SOS hierarchy (\ref{resolution-3}) for which 
finite convergence usually occurs fast. For more details see \cite{interpolation}.

\subsection{In Control \& Optimal Control}
\label{sect-ocp}
Consider the Optimal Control Problem (OCP) associated with a controlled dynamical system:
\begin{equation}
\label{ocp}
\begin{array}{rl}
J^*=\displaystyle\inf_{\u(t)}&\displaystyle\int_0^T L(\x(t),\u(t))\,dt: \dot{\x}(t)\,=\,f(\x(t),\u(t)),\:t\in (0,T)\\
&\quad\x(t)\in \X,\:\u(t)\in \U,\: \forall t\in (0,T)\\
&\quad\x(0)=\x_0;\:\x(T)\in \X_T,
\end{array}
\end{equation}
where $L,f$ are polynomials, $\X,\X_T\subset\R^n$ and $\U\subset\R^p$
are compact basic semi-algebraic sets. In full generality
the OCP problem (\ref{ocp})
is difficult to solve, especially when state constraints $\x(t)\in \X$
are present. Given an admissible  state-control trajectory $(t,\x(t),\u(t))$,
its associated occupation measure $\phi_1$ up to time $T$ (resp. $\phi_2$ at time $T$) are defined by:
\[\phi_1(A\times B\times C)\,:=\,\int_{[0,T]\cap C} 1_{(A,B)}((\x(t),\u(t))\,dt;\quad\phi_2(D)=1_D(\x(T)),\]
for all $A\in\mathcal{B}(\X)$, $B\in\mathcal{B}(\U)$, $C\in\mathcal{B}([0,T])$, $D\in\mathcal{B}(\X_T)$. Then for every differentiable function $h:\X\times [0,T]\to\R$
\[h(T,\x(T))-h(0,x_0)\,=\,\int_0^T(\frac{\partial h(\x(t),\u(t))}{\partial t}
+\frac{\partial h(\x(t),\u(t))}{\partial \x}f(\x(t),\u(t)))\,dt,\]
or, equivalently, with $\bS:=[0,T]\times \X\times \U$:
\[\int_{\X_T} h(T,\x)\,d\phi_2(\x)\,=\,h(0,\x_0)+
\int_{\bS}(\frac{\partial h(\x,\u)}{\partial t}
+\frac{\partial h(\x,\u)}{\partial \x}f(\x,\u))\,d\phi_1(t,\x,\u).\]
Then {\em the weak formulation} of the OCP (\ref{ocp}) is the infinite-dimensional linear program:
\begin{equation}
\label{weak-ocp}
\begin{array}{rl}
\rho=\displaystyle\inf_{\phi_1,\phi_2\geq0}&\{\,\displaystyle\int_{\bS} L(\x,\u)\,d\phi_1:\\
\mbox{s.t.}&\displaystyle\int_{\X_T} h(T,\cdot)\,d\phi_2
-\displaystyle\int_\bS(\frac{\partial h}{\partial t}
+\frac{\partial h}{\partial \x}f)\,d\phi_1\,=h(0,\x_0)\\
&\forall h\in\R[t,\x]\,\}.
\end{array}
\end{equation}
It turns out that under some conditions the optimal values of (\ref{ocp}) and (\ref{weak-ocp}) are equal, i.e., $J^*=\rho$. 
Next, if one replaces ``{\em for all $h\in\R[t,\x,\u]$}" with ``{\em for all
$t^k\x^\alpha\u^\beta$}", $(t,\alpha,\beta)\in\N^{1+n+p}$", then
(\ref{weak-ocp}) is an instance of the GPM (\ref{frame-2}).
Therefore one may apply the Moment-SOS hierarchy (\ref{frame-2-sdp}). Under the conditions of Theorem \ref{th-general} one obtains the asymptotic convergence 
$\rho_d\to\rho=J^*$ as $d\to\infty$. For more details 
see \cite{ocp} and the many references therein. 

\subsection*{Robust control}
In some applications (e.g. in robust control) one is often interested in optimizing over sets of the form: 
\[\G\,:=\,\{\x\in\om_1:\: f(\x,\u)\,\geq\,0,\: \forall \u\in\om_2\},\]
where $\om_2\subset\R^p$, and $\om_1\subset\R^n$  is a simple set, in fact a compact
set such that one knows all moments of the Lebesgue measure $\lambda$ on $\om_1$.

The set $\G$ is difficult to handle because of the universal quantifier. Therefore
one is often satisfied with an inner approximation $\G_d\subset \G$, and if possible,
with (i) a simple form and (ii) some theoretical approximation guarantees.
We propose to approximate $\G$ from inside by sets of (simple) form
$\G_d=\{\x\in\om_1:p_d(\x)\geq0\}$ where $p_d\in\R[\x]_{2d}$.

To obtain such an inner approximation $\G_d\subset\G$, define $F:\om_1\to\R$, $\x\mapsto F(\x):=\displaystyle\min_\u \{f(\x,\u):\u\in\om_2\}$. Then with $d\in\N$, fixed, 
solve:
\begin{equation}
\label{def-pb-robust}
\displaystyle\inf_{p\in\R[\x]_{2d}}\int_{\om_1} (F-p)\,d\lambda:
f(\x,\u)-p(\x)\,\geq\,0,\:\forall (\x,\u)\in \om_1\times\om_2\}.
\end{equation}
Any feasible solution $p_d$ of (\ref{def-pb-robust}) is such that
$\G_d=\{\x:p_d(\x)\geq0\}\subset \G$. In (\ref{def-pb-robust}) 
$\int_{\om_1}(F-p)\,d\lambda=\Vert F-p\Vert_1$ (with $\Vert\cdot\Vert_1$ being the $L_1(\om_1)$-norm), and
\[\inf_p\int_{\om_1}(F-p)\,d\lambda\,=\,\underbrace{\int_{\om_1}F\,d\lambda}_{={\rm cte}}+
\inf_p\int_{\om_1}-p\,d\lambda\,=\,{\rm cte}-\sup_p\int_{\om_1}p\,d\lambda\,\]
and so in (\ref{def-pb-robust}) it is equivalent to maximize $\int_{\om_1}pd\lambda$.
Again the Moment-SOS hierarchy can be applied. This time one replaces the difficult positivity constraint
$f(\x,\u)-p(\x)\geq0$ for all $(\x,\u)\in \om_1\times\om_2$ with a certificate of positivity, with a degree bound on the SOS weights. That is, if
$\om_1=\{\x:g_{1,\ell}(\x)\geq0,\:\ell=1,\ldots,m_1\}$
and $\om_2=\{\u:g_{2,\ell}(\u)\geq0,\:\ell=1,\ldots,m_2\}$, then with
$d_{i,\ell}:=\lceil ({\rm deg}(\sigma_{i,\ell})/2\rceil$, one solves
\begin{equation}
\label{def-pb-robust-sdp}
\begin{array}{rl}
\rho_d=\displaystyle\sup_{p\in\R[\x]_{2d}}&\displaystyle\int_{\om_1} p\,d\lambda:
f(\x,\u)-p(\x)\,=\,\sigma_0(\x,\u)\\
&+\displaystyle\sum_{\ell=1}^{m_1}\sigma_{1,\ell}(\x,\u)\,g_{i,\ell}(\x)
+\displaystyle\sum_{\ell=1}^{m_2}\sigma_{2,\ell}(\x,\u)\,g_{i,\ell}(\u)\\
&\sigma_{i,\ell}\in\Sigma[\x,\u]_{d-d_{i,\ell}},\:\ell=1,\ldots,m_i,\:i=1,2.
\end{array}
\end{equation}
\begin{thm}[\cite{tractable}]
\label{th-robust}
Assume that $\om_1\times\om_2$ is compact and its associated quadratic module is Archimedean. Let $p_d$ be an optimal solution of (\ref{def-pb-robust-sdp}). If $\lambda(\{\x\in\om_1: F(\x)=0\})=0$ then 
$\displaystyle\lim_{d\to\infty}\Vert F-p_d\Vert_1=0$ and $\displaystyle\lim_{d\to\infty}\lambda (\G\setminus\G_d)=0$.
\end{thm}
Therefore one obtains a nested sequence of inner approximations $(\G_d)_{d\in\N}\subset\G$, 
with the desirable property that $\lambda(\G\setminus\G_d)$ vanishes as $d$ increases. For more details the interested reader is referred to \cite{tractable}.
\begin{example}
{\rm 
In some robust control problems one would like to approximate as closely as desired a 
non-convex set $\G=\{\x\in\om_1: \lambda_{\min}(\A(\x))\succeq0\}$ for some 
real symmetric $r\times r$ matrix-polynomial $\A(\x)$,
 and where $\x\mapsto \lambda_{\min}(\A(\x))$ denotes
its smallest eigenvalue. If one rewrites
\[\G=\{\x\in\om_1: \u^T\A(\x)\u\geq0,\:\forall \u\in\om_2\};\quad\om_2=\,
\{\u\in\R^r:\Vert\u\Vert=1\},\]
one is faced with the problem we have just described. In applying the above methodology
the polynomial $p_d$ in Theorem \ref{th-robust} approximates $\lambda_{\min}(\A(\x))$ 
from below in $\om_1$, and $\Vert p_d(\cdot)-\lambda_{\min}(\A(\cdot))\Vert_1\to 0$ as $d$ increases. For more details see \cite{output}.
}\end{example}

There are many other applications of the Moment-SOS hierarchy in Control, e.g. 
in Systems Identification \cite{cerone,sznaier}, Robotics \cite{tedrake}, 
for computing Lyapunov functions \cite{parrilo-2}, largest regions of attraction \cite{henrion-ieee},  to cite a few.

\subsection{Some inverse optimization problems}
In particular:
\subsubsection*{Inverse Polynomial Optimization} Here we are given a polynomial optimization 
problem $\P:\,f^*=\min\{f(\x):\x\in\om\}$ with $f\in\R[\x]_d$, and we are interested in the following issue:
Let $\y\in\om$ be given, e.g. $\y$ is the current iterate of a local minimization 
algorithm applied to $\P$. Find
\begin{equation}
\label{inverse-opt}
g^*\,=\,\displaystyle\arg\min_{g\in\R[\x]_d}\,\{\Vert f-g\Vert_1: \:g(\x)-g(\y)\,\geq\,0,\:\forall \x\in\om\,\},
\end{equation}
where $\Vert h\Vert_1=\sum_\alpha \vert h_\alpha\vert$ is the $\ell_1$-norm of coefficients of $h\in\R[\x]_d$. In other words, one searches for a polynomial $g^*\in\R[\x]_d$ as close as possible to $f$ and such that $\y\in\om$ is a global minimizer of $g^*$ on $\om$. Indeed if $\Vert f-g^*\Vert_1$ is small enough then $\y\in\om$ could be considered a satisfying solution of $\P$.
Therefore given a fixed small $\epsilon>0$,  the test $\Vert f-g^*\Vert_1<\epsilon$ could be a new stopping criterion for a local optimization algorithm, with a strong theoretical justification. 

Again the Moment-SOS hierarchy can be applied to solve (\ref{inverse-opt}) 
as positivity certificates are perfect tools to handle the positivity constraint ``$g(\x)-g(\y)\geq0$
for all $\x\in\om$". Namely with $\om$ as in (\ref{setk}), solve:
\begin{equation}
\label{inverse-relax}
\rho_t=\min_{g\in\R[\x]_d}\,\{\,\Vert f-g\Vert_1: 
g(\x)-g(\y)\,:=\,\sum_{j=0}^m\sigma_j(\x)\,g_j(\x),\quad\forall\x\,\},\end{equation}
where $g_0(\x)=1$ for all $\x$, and $\sigma_j\in\Sigma[\x]_{t-d_j}$, $j=0,\ldots,m$. Other norms are possible
but for the sparsity inducing $\ell_1$-norm $\Vert\cdot\Vert_1$, it turns out that an optimal solution $g^*$ of (\ref{inverse-relax}) has a canonical simple form. For more details the interested reader is referred to \cite{inverse-opt}.

\subsubsection*{Inverse Optimal Control} With the OCP (\ref{ocp}) in \S \ref{sect-ocp}, we now consider the following issue: {\em Given a database of admissible trajectories 
$(\x(t;\x_\tau),\u(t,\x_\tau))$, $t\in[\tau,T]$, starting in initial state $\x_\tau\in\X$ at time $\tau\in [0,T]$, 
does there exist a Lagrangian $(\x,\u)\mapsto L(\x,\u)$
such that all these trajectories are optimal for the OCP problem (\ref{ocp})?} This problem has important applications, e.g.,  in Humanoid Robotics to explain human locomotion \cite{book-robot}. 

Again the Moment-SOS hierarchy can be applied because
a weak version of the Hamilton-Jacobi-Bellman (HJB) optimality conditions is the perfect tool to state whether some given trajectory is $\epsilon$-optimal for the OCP (\ref{ocp}). Indeed given $\epsilon>0$ and an admissible  trajectory $(t,\x^*(t),\u^*(t))$, let $\varphi:[0,T]\times \X\to\R$, and $L:\X\times\U\to\R$, be such that:
\begin{equation}
\label{inv1}
\varphi(T,\x)\leq0,\:\forall\,\x\in\X;\:\frac{\partial\varphi(t,\x)}{\partial t}+
\frac{\partial\varphi(t,\x)}{\partial \x}f(\x,\u)+L(\x,\u)\,\geq\,0,\end{equation}
for all $(t,\x,\u)\in [0,T]\times\X\times\U$, and: $\varphi(T,\x^*(T))>-\epsilon$, 
\begin{equation}
\label{inv2}
\frac{\partial\varphi(t,\x^*(t))}{\partial t}+
\frac{\partial\varphi(t,\x^*(t))}{\partial \x}f(\x^*(t),\u^*(t))+L(\x^*(t),\u^*(t))<\epsilon,\end{equation}
for all $t\in [0,T]$. Then the trajectory $(t,\x^*(t),\u^*(t))$ is an 
$\epsilon$-optimal solution of the 
OCP (\ref{ocp}) with $\x_0=\x^*(0)$ and Lagrangian $L$. 
Therefore to apply the Moment-SOS hierarchy:

(i) The unknown functions $\varphi$ and $L$ are approximated 
by polynomials in $\R[t,\x]_{2d}$ and $\R[\x,\u]_{2d}$, where $d$ is the parameter in the Moment-SOS hierarchy (\ref{frame-2-sdp-dual}). 

(ii) The above positivity constraint (\ref{inv1}) on $[0,T]\times\X\times\U$
is replaced with a positivity certificate with degree bound on the SOS weights.

(iii) (\ref{inv2}) is stated for every trajectory $(\x(t;\x_\tau),\u(t,\x_\tau))$, $t\in[\tau,T]$,
in the database. Using a discretization $\{t_1,\ldots,t_N\}$ of the interval $[0,T]$, the positivity constraints (\ref{inv2}) then become a set of linear constraints on the coefficients of the unknown polynomials $\varphi$ and $L$.

(iv) $\epsilon$ in (\ref{inv2}) is now taken as a variable and one minimizes a criterion of the form 
$\Vert L\Vert_1+\gamma\, \epsilon$, where $\gamma>0$ is chosen to balance between the sparsity-inducing norm $\Vert L\Vert_1$ 
of the Lagrangian and the error $\epsilon$ in the 
weak version of the optimality conditions (\ref{inv1})-(\ref{inv2}). A detailed discussion and related results can be found in \cite{inverse-ocp}.

\subsection{Optimal design in statistics}
 In designing experiments
 one models the responses $z_1,\ldots,z_N$ of a random \textit{experiment} whose inputs are represented by a vector $\t=(t_{i})\in\R^{n}$ with respect to known \textit{regression functions} $\Phi=(\varphi_{1},\ldots,\varphi_{p})$, namely:
$z_{i}=\sum_{j=1}^{p}\theta_{j}\,\varphi_{j}(t_i)+\varepsilon_{i}$, $i=1,\ldots,N$,
where $\theta_1,\ldots,\theta_p$ are unknown parameters that the experimenter wants to estimate, $\varepsilon_{i}$ is some noise and the $(t_{i})$'s are chosen by the experimenter in a \textit{design space} $\mathcal{X}\subseteq\R^{n}$. Assume that the inputs $t_i$, $i=1,\ldots,N$, are chosen within a set of distinct points $\x_1, \ldots, \x_\ell\in\mathcal{X}$, $\ell\leq\N$, and let $n_k$ denote the number of times the particular point $\x_k$ occurs among $t_1,\ldots, t_N$. A design $\xi$ is then defined by:
\begin{equation}
\label{design}
\xi\,=\,\left(\begin{array}{ccc}\x_1&\ldots&\x_\ell\\
\frac{n_1}{\N}&\ldots&\frac{n_\ell}{N}\end{array}\right).\end{equation}
The matrix  $\M(\xi):=\sum_{i=1}^\ell \frac{n_i}{N}\Phi(\x_i)\,\Phi(\x_i)^T$ is called the information matrix
of $\xi$. Optimal design is concerned with finding a set of points
in $\mathcal{X}$ that optimizes a certain statistical criterion $\phi(\M(\xi))$, which must be  real-valued, positively homogeneous, non constant, upper semi-continuous, isotonic w.r.t. Loewner ordering,
and concave. For instance
in {\em D-optimal design} one maximizes $\phi(\M(\xi)):=\log {\rm det}(\M(\xi))$ over all $\xi$ of the form (\ref{design}). This is a difficult problem and so far most methods have used a discretization of the design space $\mathcal{X}$. 

The Moment-SOS hierarchy that we describe below does not rely an any discretization and 
works for an arbitrary compact basic semi-algebraic design space $\mathcal{X}$ as defined in (\ref{setk}). Instead we look for an atomic measure on $\mathcal{X}$ (with finite support) and we proceed in two steps:\\
$\bullet$ In the first step one solves the hierarchy of convex optimization problems indexed by 
$\delta=0,1,\ldots$.
\begin{equation}
\label{step-1}
\begin{array}{rl}
\rho_\delta=\displaystyle\sup_\y&\{\log {\rm det}(\M_d(\y)):\:y_0=1\\
&\M_{d+\delta}(\y)\,\succeq\,0;\:\M_{d+\delta-d_j}(g_j\,\y)\succeq0\},
\end{array}
\end{equation}
where $d$ is fixed by the number of basis functions $\varphi_j$ considered
(here the monomials $(\x^\alpha)_{\alpha\in\N^n_d}$). (Note that (\ref{step-1}) is
not an SDP because the criterion is not linear in $\y$, but it is still a tractable convex problem.) This provides us with an optimal solution 
$\y^*(\delta)$. In practice one chooses $\delta=0$.\\
$\bullet$ In a second step we extract an atomic measure $\mu$ from the ``moments" $\y^*(\delta)$, e.g. via Nie's method \cite{nie-3} which consists of solving the SDP:
\begin{equation}
\label{step-2}
\begin{array}{rl}
\rho_r=\displaystyle\sup_\y&\{L_\y(f_r):\:y_\alpha=y^*_\alpha(\delta),\:\forall\alpha\in\N^n_{2d}\\
&\M_{d+r}(\y)\,\succeq\,0;\:\M_{d+r-d_j}(g_j\,\y)\succeq0\},
\end{array}
\end{equation}
where $f_r$ is a (randomly chosen) polynomial strictly positive on $\mathcal{X}$. If 
$(y^*_\alpha(\delta))_{\alpha\in\N^n_{2d}}$ has a representing measure then it
has an atomic representing measure, and generically the rank condition (\ref{rank}) will be 
satisfied. Extraction of atoms is obtained via a linear algebra routine. We have tested 
this two-steps method on several  non-trivial numerical experiments (in particular with highly non-convex design spaces $\mathcal{X}$) and in all cases we were able to obtain a design.
For more details the interested reader is referred to \cite{opt-design}.

\subsection*{Other applications \& extensions}
In this partial overview, by lack of space we have not described some impressive success stories of the Moment-SOS hierarchy, e.g. 
in coding \cite{bachoc}, packing problems in discrete geometry \cite{delaat-1,vallentin}. Finally,
there is also a {\em non-commutative} version \cite{pironio} of the Moment-SOS hierarchy 
based on non-commutative positivity certificates \cite{helton} and with important applications in quantum information \cite{quantum}.

\section{Conclusion}
The list of important applications of the GPM is almost endless 
and we have tried to convince the reader that the Moment-SOS hierarchy is one
promising powerful tool for solving the GPM with already some success stories. However much remains to be done 
as  its brute force application does not scale well to the problem size. One possible research direction is to exploit symmetries and/or sparsity in large scale problems. Another one is to 
determine alternative positivity certificates which are less expensive in terms of computational burden to avoid the size explosion
of SOS-based positivity certificates.


\begin{thebibliography}{las}
\bibitem{barak}
Barak, B., Steurer, D. {\em Sum-of-Squares Proofs and the Quest toward Optimal Algorithms}, Proc. ICM 2014, Seoul, Korea.
\bibitem{bachoc}
Bachoc, C., Vallentin, F. {\em New upper bounds for kissing numbers from semidefinite programming},
J. Amer. Math. Soc. {\bf 21} (2008), 909--924.
\bibitem{sznaier}
Camps, O., Sznaier, M. {\em The interplay between Big-data and sparsity in systems identification}, in Laumond J.-P., Mansard N., and Lasserre J.B. (editors) {\em Geometric and Numerical Foundations of Movements}, pp. 133--159, Springer Tracts in Advanced Robotics {\bf 117}, Springer, New York, 2017.
\bibitem{candes}
Cand\`es, E.J., Fernandez-Granda, C. {\em Towards a Mathematical Theory of 
Super-Resolution}, Comm. Pure Appl. Math. {\bf 67} (2014), 906--956.
\bibitem{cerone}
Cerone, V., Piga, D., Regruto, D. {\em Set-membership error-in-variables identification through convex relaxation techniques},
IEEE Trans. Aut. Control {\bf 57} (2012), 517--522
\bibitem{super-resolution}
De Castro, Y.,  Gamboa, F.,  Henrion, D.,  Lasserre, J.B.
{\em Exact solutions to Super Resolution on semi-algebraic domains in higher dimensions},
IEEE Trans. Info. Theory  {\bf 63} (2017), 621--630.
\bibitem{opt-design}
De Castro, Y.,  Gamboa, F.,  Henrion, D.,  Hess, R., Lasserre, J.B.
{\em Approximate Optimal Designs for Multivariate Polynomial Regression},
LAAS report No 17044. 2017, Toulouse, France.
{\tt hal-01483490v2}. To appear in Annals of Statistics.
\bibitem{upper-bound}
De Klerk, E., Lasserre, J.B., Laurent, M., Sun Zhao. (2017)
{\em Bound-Constrained Polynomial Optimization Using Only Elementary Calculations},
Math. Oper. Res. {\bf 42} (2017),  834--853.
\bibitem{delaat-1}
de Laat, D., Vallentin, F. {\em A semidefinite programming hierarchy for packing problems in discrete geometry},
Math. Program. {\bf 151} (2015), 1--25.
\bibitem{goemans}
Goemans, M.X., Williamson, D.P.  {\em Improved approximation algorithms for maximum cut and satisfiability problems using semidefinite programming},
JACM {\bf 42} (1995), 1115--1145.
\bibitem{helton}
Helton J.W., McCullough, S. {\em A Positivstellensatz for Noncommutative
Polynomials}, Trans. Amer. Math. Soc. {\bf  356} (2004), 3721--3737
\bibitem{henrion-ieee}
Henrion, D., Korda, M. {\em Convex computation of the region of attraction of polynomial control systems},
IEEE Trans. Aut. Control {\bf 59} (2014), 297--312.
\bibitem{output}
Henrion, D., Lasserre, J.B. {\em Convergent relaxations of polynomial matrix inequalities and static output feedback}, IEEE Trans. Auto. Control {\bf 51} (2006), 192--202.
\bibitem{gloptipoly}
Henrion, D., Lasserre, J.B., Lofberg, J. {\em Gloptipoly 3: moments, optimization and semidefinite programming},  Optim. Methods Softwares {\bf 24} (2009), 761--779.
\bibitem{sirev}
Henrion, D.,  Lasserre, J.B., Savorgnan, C. {\em Approximate volume and integration of basic semi-algebraic sets}, SIAM Review {\bf 51} (2009), 722--743.  
\bibitem{interpolation}
Josz, C. Lasserre, J.B., Mourrain, B. {\em Sparse polynomial interpolation: compressed sensing, super resolution, or Prony?}, LAAS Report no 17279. 2017, Toulouse, France. {\tt arXiv:1708.06187}
\bibitem{khot}
Khot, S. {\em Innapproximability of NP-complete Problems, Discrete Fourier Analysis, and Geometry}, Proc. ICM 2010,  Hyderabad, India, 2010.
\bibitem{khot2}
Khot, S. {\em Hardness of approximation}, Proc. ICM 2014,  Seoul, Korea, 2014.
\bibitem{krivine1}
Krivine, J.L. {\em Anneaux pr\'eordonn\'es}, J. Anal. Math. {\bf 12} (1964), 307--326.
\bibitem{krivine2}
Krivine, J.L. {\em Quelques propri\'et\'es des pr\'eordres
dans les anneaux commutatifs}, C.R. Acd. Sci. Paris, Ser. I {\bf 258} (1964), 
3417--3418.
\bibitem{moments}
Landau, H.J. (Editor): {\em Moments in Mathematics}, Proc. Sympos. Appl. Math. {\bf 37} (1987).
\bibitem{lass-cr}
Lasserre, J.B. {\em Optimisation globale et th\'eorie des moments}, C.R. Acad. Sci. Paris, S\'er. I {\bf 331} (2000), 929--934.
\bibitem{lass-siopt}
Lasserre, J.B. {\em Global Optimization with polynomials and the problem of moments}, SIAM J. Optim. {\bf 11}(2001), 796--817.
\bibitem{focm}
Lasserre, J.B., Laurent, M., Rostalski, P. {\em Semidefinite characterization and computation of zero-dimensional real radical ideals},  Found. Comput. Math. {\bf 8} (2008), 607--647.
\bibitem{lass-book-icp}
Lasserre, J.B. {\em Moments, Positive Polynomials and Their Applications},
Imperial College Press, London, 2009.
\bibitem{lass-book-camb}
Lasserre, J.B. {\em An Introduction to Polynomial and Semi-Algebraic Optimization},
Cambridge University Press, Cambridge, UK, 2015.
\bibitem{tractable}
Lasserre, J.B. {\em Tractable approximations of sets defined with quantifiers},  Math. Program. {\bf 151}  (2015),  507--527.
\bibitem{aam}
Lasserre, J.B. {\em Computing Gaussian \& exponential measures of semi-algebraic sets},
Adv. Appl. Math. {\bf 91} (2017),  137--163.
\bibitem{lass-mor}
Lasserre, J.B. {\em Semidefinite programming vs. LP relaxations for polynomial programming},
Math.  Oper.  Res. {\bf 27} (2002), 347--360. 
\bibitem{aap}
Lasserre, J.B. {\em Bounds on measures satisfying moment conditions},
Annals Appl. Prob.  {\bf 12} (2002),  1114--1137. 
\bibitem{lass-sparsity}
Lasserre, J.B. {\em Convergent SDP-relaxations in polynomial optimization with sparsity}, SIAM J. Optim. {\bf 17} (2006),  822--843.
\bibitem{ocp}
Lasserre, J.B., Henrion, D., Prieur, C., Tr\'elat, E. {\em Nonlinear optimal control via occupation measures and LMI-relaxations},  SIAM J. Contr. Optim. {\bf 47} (2008), 1649--1666.
 \bibitem{inverse-opt}
 Lasserre,  J.B. {\em Inverse polynomial optimization},  Math. Oper. Res. {\bf 38} (2013), 418--436.
 \bibitem{book-robot}
 Laumond,  J.-P., Mansard, N., Lasserre, J.B. (editors).  {\em Geometric and Numerical Foundations of Movements}, Springer Tracts in Advanced Robotics {\bf 117}, 
 Springer, New York, 2017.
\bibitem{laurent-1}
Laurent, M. {\em A comparison of the Sherali-Adams, Lov\'asz-Schrijver and Lasserre relaxations for 0-1 programming},
Math. Oper. Res. {\bf  28} (2003), 470--496.
\bibitem{laurent-2}
Laurent, M. {\em Sums of squares, moment matrices and optimization over polynomials},
In Putinar M. and Sullivant S. (editors), {\em Emerging Applications of Algebraic Geometry}, pp. 157--270, IMA Volumes in Mathematics and its Applications {\bf 149}, 
Springer, New York, 2009.
\bibitem{molzahn}
Molzahn, D.K., Hiskens, I.K. {\em Sparsity-Exploiting Moment-Based Relaxations of the Optimal Power Flow Problem},  IEEE Trans. Power Systems {\bf 30} (2015),
3168--3180.
\bibitem{quantum}
Navascu\'es, M., Pironio, S., A Ac\'in, A. {\em A convergent hierarchy of semidefinite programs characterizing the set of quantum correlations},
New J. Physics {\bf 10} (2008).
\bibitem{nesterov}
Nesterov, Y. {\em Squared functional systems and optimization problems},
in Frenk H., Roos K., Terlaky T., and Zhang S. (editors), {\em High Performance Optimization}, pp. 405--440, Springer, New York, 2000.
\bibitem{nie-1}
Nie, J. {\em Optimality Conditions and Finite Convergence of Lasserre's Hierarchy}, Math. Program. Ser. A {\bf 146} (2014), 97--121.
\bibitem{nie-2}
Nie, J. {\em Certifying Convergence of Lasserre's Hierarchy via Flat Truncation}, Math. Program. Ser. A, {\bf 42} (2013), 485--510.
\bibitem{nie-3}
Nie, J. {\em The $\mathcal{A}$-Truncated $K$-Moment Problem}, 
Found. Comput. Math. {\bf 14} (2014), 1243--1276.
\bibitem{parrilo-1}
Parrilo, P. {\em Structured Semidefinite Programs and Semialgebraic Geometry Methods in Robustness and Optimization}, 
PhD Thesis, California Institute of Technology, Pasadena, CA, 2000.
\bibitem{parrilo-2}
Parrilo, P. {\em Semidefinite programming relaxations for semialgebraic problems}, Math. Program. {\bf 96} (2003), 293--320.
\bibitem{inverse-ocp}
Pauwels, E., Henrion, D., Lasserre, J.B. {\em Linear Conic Optimization for Inverse Optimal Control},  SIAM J. Control  Optim. {\bf 54} (2016), 1798--1825.
\bibitem{tedrake}
Posa, M., Tobenkin, M., Tedrake, R. {\em Stability analysis and control of rigid-body systems with impacts and friction},
IEEE Trans. Aut. Control {\bf 61}, (2016), 1423--1437.
\bibitem{pironio}
Pironio, S., M Navascu\'es, M., Ac\'in, A. {\em Convergent relaxations of polynomial optimization problems with noncommuting variables},
SIAM J. Optim. {\bf 20} (2010), 2157--2180
\bibitem{putinar}
Putinar, M. {\em Positive polynomials on compact semi-algebraic sets}, Indiana Univ. Math. J. {\bf 42} (1993), 969--984
\bibitem{schmudgen}
Schm\"udgen, K. {\em The Moment problem}, Springer, 2017.
\bibitem{vallentin}
Sch\"urmann, A., Vallentin, F. {\em Computational approaches to lattice packing and covering problems},
Discrete \& Comput. Geom. {\bf 35} (2006), 73--116
\bibitem{shor}
Shor, N.Z. {\em Nondifferentiable Optimization and Polynomial Problems}, Kluwer, Dordrecht, 1988.
\bibitem{vasilescu}
Vasilescu, F.-H. {\em Spectral measures and moment problems}, in {\em Spectral Theory and Applications}, Theta Ser. Adv. Math. {\bf 2} (2003), Theta, Bucharest, pp. 173--215.
\bibitem{waki}
Waki, K., Kim, S., Kojima, M., Muramatsu, M.
{\em Sums of squares and semidefinite program relaxations for polynomial optimization problems with structured sparsity},
SIAM J. Optim. {\bf 17} (2006), 218--242
\end{thebibliography}
\end{document}